\documentclass[12pt, reqno]{amsart}
\usepackage{amsmath, amsthm, amscd, amssymb, graphicx, color}
\usepackage{lmodern}

\usepackage{enumitem}
\usepackage{hyperref}
\newtheorem{theorem}{Theorem}[section]
\newtheorem{lemma}[theorem]{Lemma}

\newtheorem{corollary}[theorem]{Corollary}
\theoremstyle{definition}

\theoremstyle{remark}

\numberwithin{equation}{section}
\newcommand{\hp}{$ \mathcal{H}_\phi^\theta$}
\newcommand{\h}{\mathcal{H}}
\newcommand{\ho}{\overline{H_0^2}}

\begin{document}
	
	\title[Dual Truncated Hankel Operators: Characterizations and Properties]{Dual Truncated Hankel Operators: Characterizations and Properties}
	
		\author[Arup Chattopadhyay and Supratim Jana]{Arup Chattopadhyay$^*$ and Supratim Jana$^{**}$}
	\maketitle
	
	\paragraph{\textbf{Abstract}}
	
	We introduce the notion of the Dual Truncated Hankel Operator (DTHO) and provide several operator equation characterizations using the dual compressed shift operator. These characterizations are similar to classical results concerning Hankel operators and align with recent findings related to Truncated Hankel Operators (THO) \cite{GM}. Additionally, our work addresses comprehensive solutions to various operator equations encountered in studying THO and the classical Hankel operator. We have also established some fundamental operator-theoretic properties of DTHO that apply to general symbols and symbols under specific conditions.

	\vspace{0.5cm}
	\vspace{0.5cm}

	\section{\textbf{Introduction and Preliminaries}}
	
Let $Hol(\mathbb{D})$ denote the space of all holomorphic functions on the open unit disk $\mathbb{D}:= \{z \in \mathbb{C}: |z| < 1\}$. The Hardy–Hilbert space over $\mathbb{D}$, denoted by $H^2(\mathbb{D})$, is defined as
$$
H^2(\mathbb{D}) := \left\{ f \in Hol(\mathbb{D}) : f(z) = \sum_{n=0}^\infty a_n z^n \ \text{with} \ \sum_{n=0}^\infty |a_n|^2 < \infty \right\}.
$$

By considering non-tangential boundary limits on $\mathbb{T}(:=\{z \in \mathbb{C}: |z| = 1\}$), the Hardy–Hilbert space $H^2(\mathbb{D})$ can be identified with a closed subspace $H^2$ of $L^2$, consisting of functions whose negative Fourier coefficients vanish. Here, $L^2$ denotes the space of all square-integrable functions on $\mathbb{T}$ with respect to the normalized Lebesgue measure.

Let $H^{\infty}$ denote the algebra of all bounded analytic functions on $\mathbb{D}$. An inner function is an element of $H^{\infty}$ whose boundary values are unimodular almost everywhere on $\mathbb{T}$. Let $S$ denotes the shift operator on $H^2$, defined by $(Sf)(z) = z f(z)$ for $z \in \mathbb{T}$. The celebrated Beurling's theorem \cite{AB} asserts that every shift-invariant subspace of $H^2$ has the form $\theta H^2$, where $\theta$ is an inner function. As an immediate consequence, the invariant subspaces of the adjoint $S^*$ (the backward shift operator) are precisely the model spaces $K_\theta$, defined by
$$
K_\theta := H^2 \ominus \theta H^2 = H^2 \cap \theta \overline{H_0^2}, \quad \text{where } H_0^2 := z H^2.
$$
Over the past century, considerable attention has been devoted to the study of Toeplitz, Hankel, and composition operators acting on the Hardy–Hilbert space $H^2$ and its subspaces. In particular, the seminal work of Brown and Halmos \cite{BH} laid the foundation for the operator-theoretic study of Toeplitz operators. A Toeplitz operator with symbol $\phi \in L^\infty$ is defined by
$$
T_\phi : H^2 \to H^2, \quad T_\phi(h) = P(\phi h),
$$
where $P$ denotes the orthogonal projection from $L^2$ onto $H^2$, and $L^\infty$ is the Banach space of essentially bounded Lebesgue measurable functions on $\mathbb{T}$.

The Hankel operator $H_\phi$ on $H^2$ is defined by
$$
H_\phi : H^2 \to H^2, \quad H_\phi(h) = P J(\phi h), \quad \phi \in L^\infty,
$$
where $J$ denotes the \emph{flip} operator on $L^2$. The operator $J$ is both self-adjoint and unitary, and its action is given by
\begin{equation}\label{2eq1}
(Jf)(z) = \breve{f}(z) = f(\bar{z}), \quad z \in \mathbb{T}.
\end{equation}
In the literature, alternative definitions of the Hankel operator can be found. One such formulation is
$$
H_\phi : H^2 \to \overline{H_0^2}, \quad H_\phi(h) = (I - P)(\phi h),
$$
and another is
\begin{equation}\label{2Flipop}
H_\phi : H^2 \to H^2, \quad H_\phi(h) = P(\mathcal{J}(\phi h)),
\end{equation}
where $\mathcal{J}$ is a unitary operator on $L^2$, defined by
$$
\mathcal{J}(f)(z) = \bar{z}\, f(\bar{z}), \quad z \in \mathbb{T},
$$
also known as the flip operator. For a detailed and comprehensive study of these operators, the reader is referred to \cite{BH, AR, NKN, JRP, SCP}.

In 2007, Sarason initiated the study of a class of model space operators known as \emph{truncated Toeplitz operators} \cite{S1}.
A truncated Toeplitz operator (TTO), denoted by $A_\phi^\theta$, is a densely defined operator on the model space $K_\theta$ given by  
$$
A_\phi^\theta(h) = P_\theta(\phi h),~
\phi \in L^2, \; h \in K_\theta \cap L^\infty,
$$
where $P_\theta : L^2 \to K_\theta$ denotes the orthogonal projection. Since then, the operator-theoretic study of truncated Toeplitz operators and related model operators has played a central role in modern operator theory. For recent progress on the algebraic properties of truncated Toeplitz operators, we refer the reader to \cite{B1, CT, C, EF, MZ1, MZ2, SN}.

At the beginning of the last decade, Gu first investigated the algebraic properties of truncated Hankel operators (THO) in a seminal preprint. More recently, this work was revised and published jointly by Gu and Ma in \cite{GM}. They introduced the truncated Hankel operator $B_\phi$ as a densely defined operator on $K_\theta$, given by  \[B_\phi^\theta(h) = P_\theta(\mathcal{J}\phi h), ~ \phi \in L^2, \; h \in K_\theta \cap L^\infty,\] where $\mathcal{J}$ is the operator defined in \eqref{2Flipop}. In addition, within the broader framework of model space operator theory, Bessonov \cite{B2} studied the \emph{small} THO. In contrast, analogous results for the \emph{big} THO were later established by Ma, Yan and Zheng in \cite{MZ3}. 

Now consider the decompositions  
$L^2 = \overline{H_0^2} \oplus H^2 = \overline{H_0^2} \oplus K_\theta \oplus \theta H^2
~ \text{and} ~
K_\theta^\perp = \overline{H_0^2} \oplus \theta H^2.$
Let $Q_\theta$ denote the orthogonal projection from $L^2$ onto $K_\theta^\perp$. Then $Q_\theta$ can be expressed as  
\begin{equation}\label{2Qpp}
Q_\theta = I - P_\theta = (I-P) + \theta P \bar{\theta}.
\end{equation}  

The \emph{dual truncated Toeplitz operator} (DTTO) was introduced and studied by Ding and Sang \cite{DS}, and is defined by  
\[
D_\phi : K_\theta^\perp \to K_\theta^\perp, 
\qquad 
D_\phi(h) = Q_\theta(\phi h), 
\quad \phi \in L^\infty.
\]  
In \cite{DS}, Ding and Sang investigated several algebraic properties of DTTOs. 
Subsequently, in \cite{CG}, Gu provided further characterizations of DTTO using the compressed shift operator $U_\theta:= D_z$, studied several related operator equations on $K_\theta^\perp$, and generalized the algebraic properties obtained in \cite{DS}. In this sense, the study of DTTO by Gu runs parallel to the classical study of Toeplitz operators by Brown and Halmos \cite{BH}, as well as to Sarason's work on truncated Toeplitz operators \cite{S1}. For recent progress on both TTO and DTTO, we refer the reader to \cite{CKLP1, CKLP2, CR, CT, EF, DQS}.

In this article, we introduce the notion of the \emph{dual truncated Hankel operator} (DTHO), denoted by $\mathcal{H}_\phi^\theta$, and is defined as
\begin{equation}
    \mathcal{H}_\phi^\theta : K_\theta^\perp \to K_\theta^\perp, 
    \quad 
    \mathcal{H}_\phi^\theta(h) = Q_\theta J(\phi h) 
    = (I-P)J(\phi h) + \theta P\bar{\theta} J(\phi h),
\end{equation}
for $\phi \in L^\infty$. Here, $J$ denotes the flip operator defined in \eqref{2eq1}. 
We study various structural and algebraic properties of these operators in detail. 
For later use, we also denote $Q := I-P$, the orthogonal projection from $L^2$ onto $\overline{H_0^2}$.

A significant distinction in our study is that, while the truncated Hankel operator (THO) is compact for a broad class of symbols (despite its symbol not being uniquely determined), the dual truncated Hankel operator (DTHO) is compact if and only if its symbol is zero. Moreover, the DTHO is bounded if and only if its symbol is essentially bounded. These behaviors are quite similar to the classical Toeplitz operator, and they arise from the fact that part of the DTHO resembles the structure of a Toeplitz operator. This naturally leads to the following question:   
	\begin{equation}\label{2Q}
		{\text{\it Do DTHOs share any structural similarities with the classical Hankel operator }?}
	\end{equation}
We provide an affirmative answer to this question. In other words, we succeeded in showing that some algebraic properties of DTHOs follow the path of the classical Hankel operator under certain conditions. 

It is well known that a classical Hankel operator on $H^2$ satisfies the operator equation $S^*H_\phi=H_\phi S$, whereas Gu and Ma showed that $S_\theta^* B_\phi - B_\phi S_\theta$ is a specific rank two operator \cite{GM}. Here, in our study of DTHO, we show that the inherited operator equation $U_\theta^* \mathcal{H}_\phi^\theta = \mathcal{H}_\phi^\theta U_\theta$ results in a particular two-ranked operator. Moreover, the similar operator equations that occurred in the study of THO have been examined in our work as well, along with the solutions of the correlated operator equations. As a consequence, we obtain several characterizations of DTHOs (see Theorem~\ref{2Ch1},~\ref{2Ch2},~\ref{2Ch3},~\ref{2Ch4}).

The article is organized as follows. Section~2 establishes several operator equations satisfied by the dual truncated Hankel operator (DTHO). Section~3 is devoted to the analysis of these equations, where we derive a parametric representation of their solutions and investigate the intertwining properties of $\mathcal{H}_\phi^\theta$ with $U_\theta$ and $U_\theta^*$. In Section~4, we examine various algebraic properties of the DTHO and highlight its similarities with the classical Toeplitz and Hankel operators, including hyponormality.

\section{\textbf{Operator Equations of DTHO}}
	
	It is easy to check that  $\left\{\bar{z}^n, \theta z^m: ~n\in\mathbb{N};~ m\in \mathbb{N}\cup \{0\}\right\}$ forms an orthonormal basis of $K_\theta^\perp.$ For any function $F\in L^2,$ we denote $(F)_k$ as the $k^{th}$ Fourier coefficient, i.e, $ (F)_k= \langle F , e_k \rangle,$ where $\{ e_k : k\in \mathbb{Z} \}$ is the standard basis of $L^2.$ Let $U= U_\theta $ be the compressed shift operator, which is $D_z$, the DTTO with the symbol $z.$  We also use the notation $\mathcal{H}_\phi$ instead of $\mathcal{H}_\phi^\theta$ during the proofs throughout the article. For a Hilbert space $\mathcal{H}$, and for $f,g\in \mathcal{H}$, $f\otimes g$ denotes the rank one operator on $\mathcal{H}$ defined by $f\otimes g(h)=\langle h, g\rangle f$. Moreover, given an inner function $\theta$, we denote $\theta_0=(\theta)_0$, the $0^{th}$ Fourier coefficient of $\theta$, and $\phi^*(z)= \overline{\phi(\bar{z})}$, for $\phi \in L^2$.  
    
    Due to the usefulness of the action of $U_\theta$ and $U_\theta^*$ on the elements of $K_\theta^\perp$ throughout this article, it is important to recall their explicit structure. This was established by Gu (see Lemma~2.1 of \cite{CG}) and is stated as follows.

	\begin{lemma}\label{2L1}
		For $h \in K_\theta ^ \perp$, we have 
        \begin{enumerate}[label=(\roman*)]
            \item $U_\theta h =  Uh=zh+\langle h, \bar{z} \rangle(\bar{\theta}_0\theta-1),$

            \item $U_\theta^* h = U^*h= \bar{z}h +\langle h, \theta \rangle (\theta_0 - \theta)\bar{z}, ~~\text{and}$

            \item $I- U_\theta^*U_\theta= (1-|\theta_0|^2)\bar{z}\otimes \bar{z}, ~~I- U_\theta U_\theta^*=(1-|\theta _0|^2) ~\theta \otimes \theta.$
        \end{enumerate}
    \end{lemma}


Note that $U_\theta$ acts as an isometry from $K_\theta^\perp \ominus \big(\vee\{\bar{z}\}\big)$ onto $K_\theta^\perp \ominus \big(\vee\{\theta\}\big)$, while $U_\theta^*$ serves as an isometry in the reverse direction, namely from $K_\theta^\perp \ominus \big(\vee\{\theta\}\big)$ onto $K_\theta^\perp \ominus \big(\vee\{\bar{z}\}\big)$. In particular, $U_\theta$ maps $\bar{z}$ to $\bar{\theta}_0 \theta$, while $U_\theta^*$ maps $\theta$ to $\theta_0 \bar{z}$. Now, we begin by computing $\mathcal{H}_\phi^\theta U_\theta - U_\theta^* \mathcal{H}_\phi^\theta$, in analogy with the corresponding investigations for the classical Hankel operator and the truncated Hankel operator. This leads us to the following conclusions.

	\begin{theorem}\label{2Th1}
		The dual truncated Hankel operator $\mathcal{H}_\phi^\theta$ satisfies the following operator equation: 
		\begin{equation}\label{2eq50}
			\mathcal{H}_\phi^\theta U_\theta - U_\theta^* \mathcal{H}_\phi^\theta= \bar{z}\otimes Q_\theta J(\alpha_{\phi^*}) - Q_ \theta J (\alpha_\phi) \otimes \bar{z}, 
		\end{equation}
		$\text{ where }  \alpha_\phi= \phi(1-\bar{\theta}_0 \theta). $
	\end{theorem}
	
	\begin{proof}
		For $h \in K_\theta ^ \perp \ominus \left(\vee\{\bar{z}\}\right)$, by using Lemma \ref{2L1}, we have $U( h )= zh,$ and
		\begin{align*}
			& U^*\mathcal{H}_\phi(h)=\bar{z}\mathcal{H}_\phi(h)+\langle \mathcal{H}_\phi(h), \theta \rangle (\theta_0 - \theta) \bar{z} \\
			& = \bar{z} [QJ(\phi h)+ \theta P \bar{\theta}J (\phi h)]+ \langle Q_\theta J(\phi h), \theta \rangle (\theta_0 - \theta) \bar{z} \\
			& = Q\bar{z} J(\phi h)- ( J(\phi h) )_0 \bar{z} + \theta P \bar{\theta} \bar{z} J (\phi h) + (\bar{\theta}J(\phi h))_0\bar{z} \theta +   (\bar{\theta}J(\phi h))_0 (\theta_0 - \theta) \bar{z} \\
			& = Q_\theta \bar{z}J(\phi h) -\{ ( J(\phi h) )_0 - \theta_0 ( \bar{\theta}J(\phi h) )_0 \}\bar{z} \\
			& =Q_\theta J(\phi zh) - \langle h,\bar{\phi}(1- \bar{\theta}_0 \breve{\theta})  \rangle \bar{z} \\
			& = \mathcal{H}_\phi (zh)- \langle h,Q_\theta J({\phi}^*(1- \bar{\theta}_0 {\theta}) ) \rangle \bar{z}.
		\end{align*}
		Thus,  $ \left(\mathcal{H}_\phi U - U^* \mathcal{H}_\phi \right)(h)= (\bar{z} \otimes Q_\theta J (\alpha_{\phi^*}))(h). $ Now,  consider $\gamma := (\mathcal{H}_\phi U - U^* \mathcal{H}_\phi) (\bar{z}) $. Therefore, for $h_1= h +c\bar{z} \in K_\theta^\perp,$ $c\in\mathbb{C}$, we have
		\begin{align*}
			&  (\mathcal{H}_\phi U - U^* \mathcal{H}_\phi) (h_1) =  (\mathcal{H}_\phi U - U^* \mathcal{H}_\phi) (h) + (\mathcal{H}_\phi U - U^* \mathcal{H}_\phi) (c\bar{z})\\
			&  = \langle h,\beta_\phi \rangle \bar{z}+ c\gamma, \text{ where } \beta_\phi = Q_\theta J (\alpha _ {\phi^*}) \\
			&  = \langle h +c\bar{z} ,\beta_\phi \rangle \bar{z} -c \langle \bar{z}, \beta_\phi \rangle \bar{z} + c\gamma \\
			&  = \langle h_1, \beta_\phi \rangle \bar{z} -\langle h_1, \bar{z}\rangle \langle \bar{z}, \beta_\phi \rangle \bar{z} + \langle h_1, \bar{z} \rangle \gamma\\
			&  = \langle h_1, \beta_\phi - \langle \beta_\phi, \bar{z} \rangle \bar{z} \rangle \bar{z} + \langle h_1, \bar{z} \rangle \gamma  \\
			&  = [\bar{z} \otimes (\beta_\phi - \langle \beta_\phi, \bar{z} \rangle \bar{z})](h_1) + (\gamma \otimes \bar{z})(h_1) .
		\end{align*}
		Next, we want to compute $ \gamma = (\mathcal{H}_\phi U - U^* \mathcal{H}_\phi )(\bar{z}).  $  Again by using Lemma \ref{2L1}, we have
		$\mathcal{H}_\phi U (\bar{z})= \mathcal{H}_\phi(\bar{\theta}_0\theta)=  Q_\theta J (\phi\bar{\theta_0} \theta  )$, and
		\begin{align*}
			& U^*\mathcal{H}_\phi(\bar{z})= \bar{z}\mathcal{H}_\phi(\bar{z})+ \langle \mathcal{H}_\phi (\bar{z}) ,\theta \rangle (\theta_0 - \theta)\bar{z} \\
			&  = \bar{z} [ QJ(\phi \bar{z}) + \theta P \bar{\theta}J (\phi \bar{z}) ] +\langle Q_\theta J(\phi \bar{z})  , \theta \rangle(\theta_0 - \theta)\bar{z} \\
			&  =  Q\bar{z}J(\phi \bar{z}) - (J(\phi \bar{z}))_0 \bar{z} + \theta P \bar{\theta}\bar{z}J (\phi \bar{z}) + (\bar{\theta}J(\phi \bar{z}))_0 \bar{z}\theta + (\bar{\theta}J(\phi \bar{z}))_0 (\theta_0 - \theta)\bar{z} \\
			&  = Q_\theta (\phi)- \{ (J(\phi \bar{z}))_0 - \theta_0 (\bar{\theta}J(\phi \bar{z}))_0 \} \bar{z} = Q_\theta (\phi) - \delta \bar{z},
		\end{align*}
		where $\delta = \{ (J(\phi \bar{z}))_0 - \theta_0 (\bar{\theta}J(\phi \bar{z}))_0 \} .$ Therefore, $\gamma= Q_\theta J (\phi(\bar{\theta}_0\theta -1)) + \delta \bar{z},$ and hence 
		$ \mathcal{H}_\phi U - U^* \mathcal{H}_\phi = (\bar{z} \otimes (\beta_\phi - \langle \beta_\phi, \bar{z} \rangle \bar{z})) + ( Q_\theta J (\phi(\bar{\theta}_0\theta -1)) + \delta \bar{z}) \otimes \bar{z}.$
		\vspace{0.1in}
		
		On the other hand, $ - \langle \bar{z}, \beta_\phi \rangle = \theta_0(\phi \theta^*)_1 - (\phi)_1 = \delta, $ and  finally we arrive at
		$$  \mathcal{H}_\phi^\theta U_\theta - U_\theta^* \mathcal{H}_\phi^\theta= \bar{z}\otimes Q_\theta J(\alpha_{\phi^*}) - Q_ \theta J (\alpha_\phi) \otimes \bar{z},$$ 
		where  $\alpha_\phi= \phi(1-\bar{\theta}_0 \theta). $  This completes the proof. 
	\end{proof}
	As a consequence of the above theorem, we conclude that the following intertwining property holds. 	
	\begin{corollary}\label{2C1}
		Let  $Z_\theta = \left\{ \frac{dz+\tilde{h}}{1-\bar{\theta}_0 \theta} : \tilde{h} \in  J(K_\theta), d \in \mathbb{C} \right\}.$ Then for $\phi, \phi^* \in Z_\theta,$ we have $U_\theta^*\mathcal{H}_\phi^\theta = \mathcal{H}_\phi^\theta U_\theta $.
	\end{corollary}
	
	\begin{proof}
		From Theorem \ref{2Th1}, we say that $\mathcal{H}_\phi U - U^* \mathcal{H}_\phi=0$ if and only if 
		$$ \bar{z}\otimes Q_\theta J(\alpha_{\phi^*}) - Q_\theta J (\alpha_{\phi}) \otimes \bar{z} = 0 \iff Q_\theta J (\alpha_{\phi}) \otimes \bar{z} = \bar{z}\otimes Q_\theta J (\alpha_{\phi^*}) .$$ 
		Therefore, there exists a complex number $c$ such that
		\begin{align*}
			& Q_\theta J (\alpha_{\phi}) = c\bar{z}, \text{ and } Q_\theta J (\alpha_{\phi^*}) = \bar{c}\bar{z}\\
			\implies & \alpha_{\phi}= {c}z + h_1(\bar{z}) \text{ and } \alpha_{\phi^*}= \bar{c}z + h_2(\bar{z}) , \text{ for some $h_1,h_2 \in K_\theta$ }\\
			\implies & \phi (1-\bar{\theta}_0 \theta) = cz + h_1({\bar{z}}) \text{ and } \phi^* (1-\bar{\theta}_0 \theta) = \bar{c}z + h_2(\bar{z})  ,
		\end{align*}
		And hence the result follows.
	\end{proof}
	
	\begin{itemize}\label{2imp1}
		\item In particular, if $\theta_0=0,$ then $ \alpha_\phi = \phi, $ and hence  $ \mathcal{H}_\phi^\theta U_\theta - U_\theta^* \mathcal{H}_\phi^\theta= \bar{z}\otimes Q_\theta (\bar{\phi}) - Q_ \theta  (\breve{\phi}) \otimes \bar{z}.$
	\end{itemize}
	
	Next, we investigate the analogous study of the operator equation for DTHO that appeared in \cite{GM}. And we have the following observation.
	\begin{theorem}\label{2Th2}
		The dual truncated Hankel operator $\mathcal{H}_\phi^\theta$ satisfies the following operator equation:
		\begin{equation}\label{2eq25}
			\mathcal{H}_\phi ^ \theta- U_\theta^* \mathcal{H}_\phi ^ \theta U_\theta^*=  \bar{z} \otimes \left( \beta_\phi- \langle \beta_\phi, \theta \rangle \theta   \right) + \left( Q_\theta J \eta_\phi + \delta \bar{z} \right) \otimes \theta,
		\end{equation}
		where $\beta_\phi = Q_\theta (\bar{\phi}z (1-\bar{\theta}_0\breve{\theta})), \eta_\phi = \phi(\theta-\theta_0)$, and $\delta = \theta_0\left\{( J(\phi \bar{z}))_0 - \theta_0 ( \bar{\theta}J(\phi \bar{z}))_0 \right\} $. 
	\end{theorem}
	
	\begin{proof}
		For $h \in K_\theta ^ \perp \ominus \left(\vee\{\theta\}\right)$, by using Lemma \ref{2L1},  we have $U^*( h )= \bar{z}h,$ and 
		\begin{align*}
			& U^*\mathcal{H}_\phi U^* (h)= U^*\mathcal{H}_\phi(\bar{z}h)= \bar{z}\mathcal{H}_\phi(\bar{z}h)+ \langle \mathcal{H}_\phi(\bar{z}h), \theta \rangle ({\theta_0}-\theta )\bar{z}\\
			& = \bar{z} \{ QJ(\bar{z}h)+ \theta P \bar{\theta}J (\bar{z}h) \} + \langle Q_\theta J(\phi\bar{z}h),\theta \rangle ({\theta_0}-\theta)\bar{z}\\
			& = Q \bar{z}J (\phi \bar{z}h) -( (J (\phi \bar{z}h))_0) \bar{z} + \theta P \bar{\theta}\bar{z} J (\phi\bar{z}h)+ ((\bar{\theta}J (\phi \bar{z}h))_0)\theta \bar{z} + (\bar{\theta}J (\phi \bar{z}h))_0) (\theta_0-\theta)\bar{z}\\
			& = Q_\theta \bar{z}J (\phi \bar{z}h) - \{(J (\phi \bar{z}h))_0 - \theta_0(\bar{\theta}J (\phi \bar{z}h))_0 \} \bar{z}\\
			& = Q_\theta J (\phi h) - \langle h, \bar{\phi}z (1-\bar{\theta}_0\breve{\theta}) \rangle \bar{z} \\
			& = \mathcal{H}_\phi (h) - \langle h, \beta_\phi \rangle \bar{z}, \text{ where }\beta_\phi = Q_\theta (\bar{\phi}z (1-\bar{\theta}_0\breve{\theta}))).
		\end{align*}
		Therefore, $(\mathcal{H}_\phi - U^* \mathcal{H}_\phi U^*)(h)= (\bar{z} \otimes \beta_\phi)(h).  $ Now, we set $\gamma = (\mathcal{H}_\phi - U^* \mathcal{H}_\phi U^*)(\theta) $. Note that for $h_1 \in K_\theta^\perp$, we have $h_1= h+ c\theta,$ for some $c\in\mathbb{C}.$
		Thus,
		\begin{align*}
			& (\mathcal{H}_\phi - U^* \mathcal{H}_\phi U^*)(h_1)= \langle h, \beta_\phi\rangle\bar{z}+ c\gamma \\
			& = \langle h+c\theta, \beta_\phi\rangle\bar{z}- c\langle \theta, \beta_\phi \rangle \bar{z} + c\gamma \\
			& = \langle h_1,( \beta_\phi - \langle \beta_\phi, \theta \rangle \theta) \rangle\bar{z} + \langle h_1, \theta \rangle \gamma\\
			& = [\bar{z} \otimes ( \beta_\phi - \langle \beta_\phi, \theta \rangle \theta)](h_1) + (\gamma \otimes \theta) (h_1).
		\end{align*}
		Next, we want to calculate $\gamma =(\mathcal{H}_\phi - U^* \mathcal{H}_\phi U^*)(\theta).$ Again by using Lemma \ref{2L1}, we have $  \mathcal{H}_\phi (\theta)= Q_\theta J(\theta\phi),$
		and 
		\begin{align*}
			& U^*\mathcal{H}_\phi U^* (\theta)= U^* \mathcal{H}_\phi({\theta_0}\bar{z}) = \bar{z} \mathcal{H}_\phi({\theta_0}\bar{z}) +\langle \mathcal{H}_\phi({\theta_0}\bar{z}), \theta \rangle(\theta_0 - \theta)\bar{z}\\
			& = \bar{z}\{ QJ(\theta_0\bar{z}\phi) + \theta P \bar{\theta}J(\theta_0\bar{z}\phi) \} + \langle Q_\theta J(\phi{\theta_0}\bar{z}), \theta \rangle(\theta_0 - \theta)\bar{z}\\
			& = Q_\theta \bar{z} J ( \theta_0\bar{z}\phi ) - \{ (J(\theta_0\bar{z}\phi))_0 - \theta_0 (\bar{\theta}J(\theta_0\bar{z}\phi))_0  \} \bar{z}\\
			& = Q_\theta J(\theta_0 \phi) - \delta \bar{z}.
		\end{align*}
		This yields $\gamma = Q_\theta J(\phi(\theta-\theta_0)) + \delta \bar{z} = Q_\theta J(\eta_\phi) + \delta \bar{z},$ and hence we have the desired equality. This completes the proof. 
	\end{proof}
	In particular, we have the following observation.	
	
	\begin{corollary}\label{2C2}
		If $\theta_0 = 0 ,$ then the operator \hp ~satisfies 
		$ \mathcal{H}_\phi ^ \theta- U_\theta^* \mathcal{H}_\phi ^ \theta U_\theta^*=  \bar{z} \otimes \left( \beta_\phi- \langle \beta_\phi, \theta \rangle \theta   \right) + \left( Q_\theta J (\phi \theta)) \right) \otimes \theta.$
	\end{corollary}
	
	\begin{proof}
		If $\theta_0=0,$ then $\delta=0,$ and $\gamma= Q_\theta J(\phi\theta) .$ Then the rest follows from the above theorem.
        \end{proof}
	In the sequel, we investigate two more operator equations concerning DTHO, as done in \cite{GM} for THO. 
	
	\begin{theorem}\label{2Th3}
		The operator \hp ~satisfies the following operator equation:
        \begin{equation}\label{2eqn1}
            \mathcal{H}_\phi ^ \theta - U\mathcal{H}_\phi ^ \theta U =\theta \otimes \left(Q_\theta J (\alpha_{\phi^*}) - \langle Q_\theta J (\alpha_{\phi^*}), \bar{z} \rangle \bar{z}\right) + \left(Q_\theta J ( \phi \bar{z}(1-\bar{\theta}_0 \theta) ) + \delta \theta\right) \otimes \bar{z},
        \end{equation} where $ \alpha_\phi = \phi(\theta - \theta_0), $ and $\delta = \bar{\theta}_0 \langle \phi \theta ( \theta^*- \bar{\theta}_0 ), z \rangle. $
		
	\end{theorem}
	
	\begin{proof}
		The proof is similar to the proof of Theorem \ref{2Th1}. First, we compute $ \left(\mathcal{H}_\phi - U\mathcal{H} U \right)(h)$ for $h\in K_\theta^\perp \ominus \left(\vee\{ \bar{z} \}\right)$, and then we calculate the action of $ (\mathcal{H}_\phi - U\mathcal{H} U)$ on $h_1,$ where $h_1= h+c\bar{z} $ for some complex number $c.$ We leave the details of the proof to the reader.
	\end{proof}
	
	\begin{theorem}\label{2Th4}
		The operator \hp ~satisfies the following operator equation:
        \begin{equation}\label{2eqn2}
            \mathcal{H}_\phi ^ \theta U^* - U \mathcal{H}_\phi ^ \theta= \theta \otimes Q_\theta J (\alpha_{\phi^*})- Q_\theta J(\alpha_ \phi) \otimes \theta, \end{equation} where $\alpha_\phi = \phi \bar{z}(\theta - \theta_0).$
	\end{theorem}
	
	\begin{proof}
		Since the proof proceeds analogously to that of Theorem \ref{2Th2}, we omit the details.
	\end{proof}

	\section{\textbf{The Solutions of Hankel Related Operator Equations}}

    This section investigates various operator equations that emerged in the study of the dual truncated Hankel operator in the previous section, with a focus on examining how closely their solutions are related to DTHOs.

  The following lemma has been used extensively throughout this section. Gu employed the identities stated below in \cite{CG} while studying operator equations related to the DTTO. We present these identities in the form of the following lemma for completeness.
	\begin{lemma} \label{2L2}
		The compressed shift operator $U_\theta$ and its adjoint $U_\theta^*$ are invertible when $\theta_0 \neq 0$ and satisfy the following identities.
    \begin{align*}
			& (i) ~U_\theta^n( \bar{z}\bar{\eta} ) = Q( \bar{z}\bar{\eta} . z^n ) + \bar{\theta}_0\theta P(\bar{z}\bar{\eta} z^n) ~(\text{ for } \eta \in H^\infty) ,\\
			& (ii) ~{U_\theta ^*}^n (\theta\psi) = \theta_0 Q(\psi. \bar{z}^n) + \theta P(\psi. \bar{z}^n)~ (\text{ for } \psi \in H^\infty),\\
			& (iii) ~U_\theta^{-n} (\theta \psi)= \frac{1}{\bar{\theta}_0} Q(\psi \bar{z}^n) + \theta P(\psi \bar{z}^n) ~(\text{ for } \psi \in H^\infty, \theta_0\neq 0  ), \\
			& (iv) ~{U_\theta^*}^{-n} (\bar{z}\bar{\eta}) = Q(\bar{z}\bar{\eta} . z^n) + \frac{\theta}{\theta_0} P (\bar{z}\bar{\eta} .z^n)~ (\text{ for } \eta \in H^\infty, \theta_0\neq 0  ), \\
			& (v) ~{U_\theta^*}^{-1} \bar{z}^{n+1}= \bar{z}^n (n \geq 1);   {U_\theta^*}^{-1}(\theta z^m) = \theta z^{m+1}(m\geq 0); {U_\theta^*}^{-1}(\bar{z}) = \theta/\theta_0~ (\theta_0\neq 0 ).
		\end{align*}
	\end{lemma}
	\begin{proof}
		One can apply the induction principle on the action of $U_\theta, U_\theta^{-1}, U_\theta^*$, and ${U_\theta^*}^{-1}$ on the elements of $K_\theta^\perp$ to get these identities. The proofs are left to the reader.
	\end{proof}
	We have already shown in Section 2 that $ \mathcal{H}_\phi^\theta-U_\theta^* \mathcal{H}_\phi^\theta U_\theta^*$ is a rank-two operator. Now the question is: For $A \in \mathcal{B}(K_{\theta}^{\perp})$, does $ A- U_\theta^* A U_\theta^*=0$ have a non-trivial solution or not?

	In 2002, Martinez-Avendano studied operator equations over the Hardy space $H^2$ in this direction. In fact, in Corollary 2.2 of \cite{RAMA}, he showed that the operator equation $ X=S^*XS^*$, where $S$ is the shift operator, has only the trivial solution. But in our case, that is, over the space $K_\theta^\perp$,  we prove the existence of non-trivial solutions as mentioned in the following theorem.
	
	\begin{theorem} \label{2Th5}
		Let $A$ be bounded operator on $K_\theta^\perp$ satisfying the operator equation $ A= U_\theta^* A U_\theta ^* $ . Then,
		\vspace{0.1in}
		
		(i)  if $\theta_0 = 0$, then there exists an $L^\infty$- function $ \psi $ such that $A(h)= P\mathcal{H}_{\psi \breve{\theta}} (Qh) $,
		\vspace{0.1in}

		(ii) if $\theta_0 \neq 0,$ then there exist  $\psi, \eta\in H^\infty$ such that $A= F_{\eta} + G_\psi$, where $F_{\eta}$ and $G_\psi$ are operators defined on $K_\theta ^\perp$ and given by 
		
		$$ F_{\eta}h = \theta_0 \bar{z}^2 \bar{\eta} \theta^* J(Ph) + Q( \bar{z}^2 \bar{\eta}J(Qh)) + \frac{\theta}{\theta_0} P ( \bar{z}^2 \bar{\eta} J(Qh)  ), $$

		$$ G_\psi h = \theta \psi \bar{z} J (Qh) + \theta_0 \left[ \theta_0 Q(\psi \bar{z} \theta^* J (Ph)) + \theta P (\psi \bar{z} \theta^* J (Ph)) \right]. $$
		
	\end{theorem}
	\vspace{0.5cm}
	
	\begin{proof}
		Case-1 ($\theta_0=0$): In this case by Lemma~\ref{2L1}, we have $U(\bar{z})=0$ and $U^*(\theta)=0.$ Therefore, $ A\theta = (U^*AU^*)\theta=U^*A 0=0,$ which implies $Az^n\theta= {U^*}^{n} A {U^*}^{n}(z^n \theta)= {U^*}^n A \theta=0 $, and hence $A|_{\theta H^2} = 0.$
		\vspace{0.1in}

		Again by using Lemma~\ref{2L1}, we conclude $ \langle A\bar{z}^n, \bar{z}^m \rangle = \langle {U^*}^m A {U^*}^m \bar{z}^n , \bar{z}^m \rangle = \langle A\bar{z}^{m+n} , U^m{\bar{z}^m} \rangle = \langle A\bar{z}^{m+n} , 0 \rangle = 0.$ 
		\vspace{0.1 in}  
		
		Therefore,  $A$ can be considered as an operator acting on $\overline{H_0^2}$ and having the range in $\theta H^2.$
		\vspace{0.1 in} 
		
		For $n\geq1$, let $A\bar{z}^n = \theta h_n, $ for some $h_n \in H^2.$
		\vspace{0.1 in}
		
		So, we have $\theta h_n = A\bar{z}^n =  U^*AU^* \bar{z}^n= U^*A\bar{z}^{n+1} = U^* (\theta h_{n+1})= Q(\bar{z}\theta h_{n+1}) + \theta P \bar{\theta}(\bar{z} \theta h_{n+1})$$= 0 + \theta P (\bar{z}h_{n+1}) = \theta P (\bar{z}h_{n+1})$.
		\vspace{0.1 in} 
		
		Now, we define an operator $B: \overline{H_0^2} \rightarrow H^2  $ by $B\bar{z}^n = h_n= P(\bar{z}h_{n+1}) =  P(\bar{z}B \bar{z}^{n+1}) = T_z^* B \mathbb{S} (\bar{z}^n) $ (where $\mathbb{S}(h) = \bar{z} h$, for $h\in \ho$). Therefore, we conclude $B=B_\psi$ for some $\psi \in L^\infty$, where $B_\psi: \ho \rightarrow H^2$ is such that \begin{equation}\label{Toep}
            B_\psi(h) = PJ(\psi f).
        \end{equation}
		
		Therefore, $ A \bar{z}^n = \theta B \bar{z}^n = \theta B_\psi \bar{z}^n = \theta P J (\psi \bar{z}^n) \implies A(\bar{z}\bar{f}) = \theta P J (\psi \bar{z}\bar{f})) $ for $f\in H^2,$ and hence $A(h)= A(Qh) + A(Ph) = \theta P \bar{\theta} J(\psi \breve{\theta} Qh )+ 0=P\mathcal{H}_{\psi\breve{\theta}} (Qh) $.
		\vspace{0.1 in}
		
		Case-2 ($\theta_0 \neq 0$): Again  by using Lemma~\ref{2L1}, we have $U\bar{z}= \bar{\theta_0} \theta $, and $U^* \theta = \theta_0\bar{z}.$ Since, $A$ is a self map on $K_\theta^\perp $, we write $A\bar{z}= \overline{z \eta} +\theta \psi$ for some $\eta,\psi \in H^2.$ 
		\vspace{0.1 in}
		
		Now, we assume that $ A\bar{z}= \overline {z\eta}  \implies A\bar{z}=  U^*AU^*\bar{z} =\cdots = {U^*}^n A {U^*}^n\bar{z} = \overline {z\eta}$,
		which further by using Lemma~\ref{2L2}  implies
        \begin{align*}
             A\bar{z}^{n+1} & =  {U^*}^{-n} (\overline {z\eta}) \\
             & = Q(\overline{z\eta}. z^n) + \frac{\theta}{\theta_0} P (\overline{z\eta}.z^n)\\
             & = Q(\bar{z}^2 \eta J(\bar{z}^{n+1}) ) + \frac{\theta}{\theta_0} P (\bar{z}^2 \bar{\eta}J(\bar{z}^{n+1})) 
        \end{align*}

		
		Similarly, by using Lemma~\ref{2L1}, we have $A\theta = (U^*AU^*)\theta = U^*A(\theta_0\bar{z})= \theta_0 U^*(\overline{zw}) = \theta_0\overline{zw}. \bar{z},$ and $ A(z\theta)= (U^*AU^*) z\theta = U^* A \theta= \theta_0 \overline{zw}. \bar{z}^2,\ldots, A(z^n\theta) = \theta_0\overline{zw}. \bar{z}^{n+1}= \theta_0\overline{zw}.\bar{z}\theta^* J({\theta z}^n ). $ 
		Therefore, 
		
		\begin{equation}\label{2eq20}
			A(h) = AQh + APh = F_{\eta} (h).
		\end{equation}
		
		Next, we assume $A\bar{z}=\theta \psi\implies  A\bar{z} = U^*AU^*\bar{z}=\cdots= {U^*}^n A{U^*}^n \bar{z} = \theta \psi$, 
		which by using Lemma~\ref{2L2} implies
        \begin{align*}
            A(\bar{z}\bar{z}^n) & = {U^*}^{-n}(\theta \psi)\\
            & = \theta \psi z^n\\
            & = \theta \psi J(\bar{z}^n)\\
            & = \theta \psi \bar{z} J(\bar{z}^{n+1})
        \end{align*} and it shows $\psi \in H^{\infty}$.
        
		
		Again by using Lemma~\ref{2L2}, we get
        \begin{align*}
            A z^n\theta & = \theta_0 {U^*}^{n+1} (\theta \psi)\\
            & = \theta_0 [\theta_0 Q(\psi \bar{z}^{n+1}) + \theta P (\psi \bar{z}^{n+1})]\\
            & = \theta_0[ \theta_0 Q (\psi \bar{z}\theta^* J( z^n \theta ) ) + \theta P ( \psi \bar{z} \theta^* J(z^n \theta) )]
        \end{align*}

		Therefore,
		
		\begin{equation}\label{2eq21}
			A(h) = AQh + APh = G_\psi (h).
		\end{equation}
		Finally, combining \eqref{2eq20} and \eqref{2eq21}, we conclude $Ah = F_w h + G_\psi h.$ This completes the proof. 
	\end{proof}
	
	In particular, we have the following observation. 
	
	\begin{corollary}\label{2C5}
		Suppose $A$ is a bounded linear operator on $K_\theta ^ \perp.$ Then, $A=0$ if and only if $ A= U_\theta^* A U_\theta^*  $, $A\bar{z}=0,$ and $ A^*\theta=0.$    
		
	\end{corollary}
	
	\begin{proof}
		If $\theta_0\neq 0,$ then by the calculation done in the above theorem, it is very clear that $U^* A U^* = A$ and $A\bar{z}= 0$ implies $\psi=0=\eta,$ and hence $ A=0. $
		
		Now,  if $\theta_0=0,$ then $U^*AU^*=A ,$ and $A\bar{z} = 0= A^*\theta,$ which implies that $ P( \breve{\psi} z )=0 $ and $ Q( \bar{\psi} JP( \bar{\theta} \theta ) )=0,$ and it leads to the fact that $ \psi \in z^2 H^2 \cap \overline{H^2} \implies \psi=0$. Therefore, $A$ is the zero operator. 
	\end{proof}
	The following theorem \emph{characterizes} DTHO in several cases.
	
	\begin{theorem}\label{2Ch1}
		Let $A \in \mathcal{B}(K_{\theta}^{\perp})$ satisfying the equation $ A- U_\theta^* A U_\theta^*=  \bar{z} \otimes \left( \beta_\phi- \langle \beta_\phi, \theta \rangle \theta   \right) + \left( Q_\theta J \eta_\phi + \delta \bar{z} \right) \otimes \theta$, where $\beta_{\phi},$ $\eta_{\phi}$, and $\delta$ are as in equation  \eqref{2eq25}. Then,
		\begin{itemize}
			\item if $\theta_0\neq 0$, then $A=\mathcal{H}_{\phi}^{\theta}$ if and only if  $A\bar{z}=\mathcal{H}_{\phi}^{\theta}\bar{z}$, and
			\item  if $\theta_0=0$, then $A=\mathcal{H}_{\phi}^{\theta}$ if and only if $A\bar{z}=\mathcal{H}_{\phi}^{\theta}\bar{z}$ and $A^*\theta=(\mathcal{H}_{\phi}^{\theta})^*\theta$.
		\end{itemize}
	\end{theorem}
	\begin{proof}
		The proof follows from Theorem~\ref{2Th2} and Corollary~\ref{2C5}. 
	\end{proof}
	
	In the same spirit, Martinez-Avendano in (\cite{RAMA}, Corollary 2.4) also showed that the zero operator is the only solution to the operator equation $ X=SXS$. In a manner similar to the above, we obtain non-trivial solutions to the associated operator equation over the space $K_{\theta}^\perp$; these are presented as follows.

	\begin{theorem}
		Let $A$ be a bounded operator on $K_\theta^\perp$  satisfying the operator equation $ A= U_\theta A U_\theta $.  Then,
		\vspace{0.1in}
		
		(i) if $\theta_0 = 0$, then there exists an $L^\infty$ function $ \psi $ such that $A(h)= QJ(\phi Ph) = Q \mathcal{H}_\psi ^ \theta (Ph) $,
		\vspace{0.1in}
		
		(ii) if $\theta_0 \neq 0,$ then there exist $\psi, \eta \in H^\infty$ such that $A= F_\eta + G_\psi$ where $F_\eta$ and $G_\psi$ are operators defined on $K_\theta ^\perp$ and given by 
		
		$$ F_{\eta}h = \theta^* \bar{z}\bar{\eta}J(Ph) + \bar{\theta}_0 [ Q( \bar{z}\bar{\eta} J(Qh) + \bar{\theta}_0 \theta P(\bar{z}\bar{\eta} J(Qh) ) ], $$

		$$ G_\psi h = \bar{\theta}_0 \theta \psi J(Qh) + \frac{1}{\bar{\theta}_0} Q(\psi \theta^* J (Ph)) + \theta P (\psi \theta^* J (Ph)) . $$
	\end{theorem}
	
	\begin{proof}
		The proof of the theorem is similar to the proof of Theorem \ref{2Th5}. Thus, the details are excluded.
	\end{proof}

	\begin{corollary}
		Suppose $A$ is a bounded linear operator on $ K_\theta^\perp.$ Then, $A=0$ if and only if $ A= U_\theta A U_\theta  $, $A^*\bar{z}=0,$ and $A\theta=0.$
	\end{corollary}

    The following is another characterization of DTHO.

    \begin{theorem}\label{2Ch2}
        Let $A \in \mathcal{B}(K_{\theta}^{\perp})$ satisfying the equation $ A- U_\theta A U_\theta= \theta \otimes \left(Q_\theta J (\alpha_{\phi^*}) - \langle Q_\theta J (\alpha_{\phi^*}), \bar{z} \rangle \bar{z}\right) + \left(Q_\theta J ( \phi \bar{z}(1-\bar{\theta}_0 \theta) ) + \delta \theta\right) \otimes \bar{z}, $ where $\alpha_\phi$ and $\delta$ are as in equation \eqref{2eqn1} Then,
		\begin{itemize}
			\item if $\theta_0\neq 0$, then $A=\mathcal{H}_{\phi}^{\theta}$ if and only if  $A\theta=\mathcal{H}_{\phi}^{\theta}\theta$, and
			\item  if $\theta_0=0$, then $A=\mathcal{H}_{\phi}^{\theta}$ if and only if $A\theta=\mathcal{H}_{\phi}^{\theta}\theta$ and $A^*\bar{z}=(\mathcal{H}_{\phi}^{\theta})^*\bar{z}$.
		\end{itemize}
    \end{theorem}
    
    Let $I(U_\theta, U_\theta^*)$ denote the collection of all bounded linear operators on $K_\theta^\perp$ that intertwine $U_\theta$ and $U_\theta^*$; that is, if $X \in I(U_\theta, U_\theta^*)$, then $U_\theta^* X = X U_\theta.$ In the classical Hardy space $H^2$, the solutions to the intertwining relation $S^* X = X S$ are precisely the Hankel operators. For the model space $K_\theta$, the corresponding intertwining class has been studied by Gu and Ma in \cite{GM} and by B, Bala, Panja, and Sarkar in \cite{Jayetall}. The following theorem extends this discussion to the context of $K_\theta^\perp$.
    
	\begin{theorem}\label{2Th6}
		Let $A$ be a bounded operator on $K_\theta ^\perp$ satisfying the operator equation $ U_\theta^* A= A U_\theta $.  Then,
		\vspace{0.1in}
		
		(i) if $\theta_0 = 0$, then there exist $\Phi,\psi,\eta \in H^\infty$ such that 
		$$ Ah= A(Qh)+(APh)=  \theta PJ( \Phi z Qh ) + \overline{z\eta} \theta^* J(Ph)  + \theta P (\psi \theta^* J(Ph)),$$
		\vspace{0.1in}
		
		(ii) if $\theta_0 \neq 0,$ then there exist $\psi, \eta\in H^\infty$ such that
		
		$$ Ah= A(Ph)+A(Qh) = \overline{z\eta} \theta^* J(Ph) + \theta_0 Q(\psi \theta^* J(Ph)) + \theta P (\psi \theta^* J(Ph) $$ $$ \hspace{1.7in} + \psi \theta J(Qh) + \bar{\theta}_0 [ Q(\overline{z\eta} J (Qh)) + \frac{\theta}{\theta_0} P (\overline{z\eta} J (Qh)) ] .$$
		
	\end{theorem}
	
	\begin{proof}
		The proof can be done similarly to Theorem \ref{2Th5}; that is, by looking at the action of $A$ on the orthonormal basis. Irrespective of the choice of the value $\theta_0$, let $A\theta = \overline{z\eta} + \theta\psi.$ Then using $ U^*A= AU $, we have the following. 
		\begin{align*}
			& A\theta = \overline{z\eta} + \theta\psi \\
			& \implies U^* A\theta=AU\theta=Az\theta = U^*(\overline{z\eta} + \theta\psi) \\
			& \implies {U^*}^nA\theta = A z^n \theta = {U^*}^n( \overline{z\eta} + \theta\psi ) = \overline{z\eta}. \bar{z}^n + \theta_0 Q(\psi \bar{z}^n ) + \theta P(\psi \bar{z}^n)\\
			& \implies Az^n\theta = \overline{z\eta}. \theta^* J(z^n\theta) + \theta_0 Q(\psi \theta^* J(z^n\theta)) + \theta P(\psi \theta^* J(z^n\theta))\\
			& \implies A Ph = \overline{z\eta}. \theta^* J(Ph) + \theta_0 Q(\psi \theta^* J(Ph)) + \theta P(\psi \theta^* J(Ph)).
		\end{align*}
		Now, if $\theta_0 \neq 0,$ then $ U,U^* $ are invertible and so $ A= {U^*}^{-1}AU ,$  and therefore we have
		\begin{align*}
			& A\bar{z}= {U^*}^{-1}AU(\bar{z})= \bar{\theta_0}{U^*}^{-1}A(\theta)\\
			& \implies A\bar{z}^2= {U^*}^{-1}AU(\bar{z}^2) = {U^*}^{-1}A(\bar{z}) = \bar{\theta_0}{U^*}^{-2}A(  \theta )\\
			& \hspace{0.2in} \vdots \\
			& \implies  A\bar{z}^n= \bar{\theta_0}{U^*}^{-n}A(\theta) = \bar{\theta}_0 [ Q(\overline{zw} J(\bar{z}^n )) + \frac{\theta}{\theta_0} P(\overline{zw} J(\bar{z}^n )) + \psi \theta J(\bar{z}^n) ].
		\end{align*}
		Finally, we conclude $$ AQh = \bar{\theta}_0 [ Q(\overline{zw} J(Qh)) + \frac{\theta}{\theta_0} P(\overline{zw} J(Qh )) + \psi \theta J(Qh) ] .$$
        
		Next, if $\theta_0=0$, then we proceed as follows:\\
        Let $$ A\bar{z} = \overline{z\xi} + \theta \phi  \quad \text{ for some }~ \xi, \phi \in H^2.$$
        Therefore, using Lemma~\ref{2L1} and ${U^*} A=AU$, we get ${U^*}^nA(\bar{z}) = {U^*}^{n-1}AU(\bar{z})= {U^*}^{n-1}A (0)= 0. $ Moreover, 
		$A\bar{z} = \overline{z\xi} + \theta \phi \implies {U^*}^n A\bar{z}= \overline{z\xi} \bar{z}^n + \theta P (\bar{z}^n \phi)=0 \implies \xi = 0,$ and $ \phi= \text{ constant }.$ Therefore, $ A\bar{z} = \theta \phi,$ and hence it is clear that the action of $A$ on ${\overline{H_0^2}} $ is contained in $\theta H^2.$
		
		Define a unitary operator $\mathcal{J} : \overline{H_0^2} \rightarrow \theta H^2$ given by $$\mathcal{J} (\bar{z}^{n+1}) = \theta z^n \quad \text{ for } n\geq 0.$$
        Then, $$ U|_{\overline{H_0^2}} = \mathcal{J}^{-1} \tau_{\bar{z}} \mathcal{J} \quad \text{ and } \quad U^*|_{\theta H^2} = \tau_{\bar{z}},$$ where $ \tau_{\phi}: \theta H^2 \rightarrow \theta H^2 $ is given by $$ \tau_{\phi}(f)= P_{\theta H^2} ( \phi f ) \quad \text{ for } \quad \phi \in L^\infty. $$ 
        Therefore, $U^*A= AU $ on $\overline{H_0^2} $ can be replaced by $\tau_{\bar{z}} A = A \mathcal{J}^{-1} \tau_{\bar{z}} \mathcal{J} \implies \tau_{\bar{z}} A \mathcal{J}^{-1} = A \mathcal{J}^{-1} \tau_{\bar{z}}.$ Thus, $ A \mathcal{J}^{-1} $ commutes with $\tau_{\bar{z}}$, and hence there exists a $ \Phi \in H^\infty $ such that $ A \mathcal{J}^{-1} = \tau_{\Phi^*}^*$. Finally, a direct calculation leads to $ AQh = \theta PJ( \Phi z Qh ).$ This completes the proof. 
	\end{proof}
	The following result is an immediate consequence of the above theorem.  
	
	\begin{corollary}\label{2C6}
		Suppose $A$ is a bounded linear operator on $K_\theta ^ \perp.$ Then, $A=0$ if and only if $ U_\theta^* A=  A U_\theta  $ and $A\theta=0=A^*\theta.$
	\end{corollary}
	
	\begin{proof}
		If $\theta_0\neq 0,$ then by the calculation done in the above theorem, it is clear that $U^* A = AU$ and $A\theta= 0$, which implies $\psi=0=\eta,$ and hence $ A=0. $
		
		Now,  if $\theta_0=0,$ then $U^*A=AU $ and $A\theta = 0= A^*\theta$  implying the fact that $ \Phi \in H^\infty \cap \overline{H_0^2} \implies \Phi = 0.$ 
	\end{proof}
	
	The following theorem characterizes DTHO in several cases.
	
	\begin{theorem}\label{2Ch3}
		Let $A \in \mathcal{B}(K_{\theta}^{\perp})$ satisfying the equation $ U_\theta^*A-A U_\theta=  \bar{z}\otimes Q_\theta J(\alpha_{\phi^*}) - Q_ \theta J (\alpha_\phi) \otimes \bar{z}$, where $\alpha_{\phi}$ as in equation  \eqref{2eq50}. Then,
		\begin{itemize}
			\item if $\theta_0\neq 0$, then $A=\mathcal{H}_{\phi}^{\theta}$ if and only if  $A\theta=\mathcal{H}_{\phi}^{\theta}(\theta)$, and
			\item  if $\theta_0=0$, then $A=\mathcal{H}_{\phi}^{\theta}$ if and only if $A\theta=\mathcal{H}_{\phi}^{\theta}(\theta)$ and $A^*\theta=(\mathcal{H}_{\phi}^{\theta})^*\theta$.
		\end{itemize}
	\end{theorem}
	\begin{proof}
		The proof follows from Theorem~\ref{2Th1} and Corollary~\ref{2C6}. 
	\end{proof}

	In the following theorem, we consider the complementary intertwining property, namely, the characterization of the solutions in $ I( U_\theta^*, U_\theta ) .$ 
	
	\begin{theorem}
		Let $A$ be bounded operator  on $K_\theta^\perp$ satisfying the operator equation $ U_\theta A= A U_\theta^*  $ Then,
		\vspace{0.1in}

		(i) if $\theta_0 = 0$, then there exist $\phi, \psi,\eta \in H^\infty$ such that 
		
		$$ Ah = A(Qh)+A(Ph )= \bar{z}\theta \psi J(Qh) + Q(\bar{z}^2 \bar{\eta} J(Qh)) + Q(\theta^* \bar{z}\phi J(Ph)), $$
		\vspace{0.01in}

		(ii) if $\theta_0 \neq 0,$ then there exist $\psi, \eta \in H^\infty$ such that
		
		$$ Ah= A(Qh) +A(Ph )  =  \bar{z}\theta \psi J(Qh) + Q(\bar{z}^2 \bar(\eta) J(Qh)) + \bar{\theta}_0 \theta P ( \bar{z}^2 \bar{\eta} J(Qh) ) $$ $$ + {\theta_0}[ \bar{z}^2 \bar{\eta} \theta^* J(Ph) + \frac{1}{\bar{\theta}_0} Q(\psi \bar{z} \theta^* J(Ph)) + \theta P (\psi \bar{z}\theta^* J(Ph)) ] .$$

	\end{theorem}

	\begin{proof}
		The proof is similar to the proof of Theorem \ref{2Th6} with appropriate changes. Therefore, the details are excluded.     
	\end{proof}

	\begin{corollary}
		Suppose $A$ is a bounded linear operator on $K_\theta ^ \perp.$ Then, $A=0$ if and only if $ U_\theta A=  A U_\theta^*  $ and $A\bar{z}=0=A^* \bar{z}.$
	\end{corollary}

 We now present an additional characterization of DTHOs.

    \begin{theorem}\label{2Ch4}
        Let $A \in \mathcal{B}(K_{\theta}^{\perp})$ satisfying the equation $ A U^* - U A= \theta \otimes Q_\theta J (\alpha_{\phi^*})- Q_\theta J(\alpha_ \phi) \otimes \theta$, where $\alpha_{\phi}$ as in equation  \eqref{2eqn2}. Then,
		\begin{itemize}
			\item if $\theta_0\neq 0$, then $A=\mathcal{H}_{\phi}^{\theta}$ if and only if  $A\bar{z}=\mathcal{H}_{\phi}^{\theta}(\bar{z})$, and
			\item  if $\theta_0=0$, then $A=\mathcal{H}_{\phi}^{\theta}$ if and only if $A\bar{z}=\mathcal{H}_{\phi}^{\theta}(\bar{z})$ and $A^*\bar{z}=(\mathcal{H}_{\phi}^{\theta})^*(\bar{z})$.
		\end{itemize}
    \end{theorem}
	\vspace{0.2cm}

	\section{\textbf{Properties of DTHO}}
	
	In this section, we apply the characterization results obtained in the previous sections for the dual truncated Hankel operator in order to investigate some of its algebraic properties. Furthermore, we demonstrate certain similarities between the DTHO, classical Toeplitz operators, and classical Hankel operators.
	
	\subsection{Toeplitzness of DTHO}

    This subsection is devoted to the study of DTHOs, focusing on the uniqueness of their symbols, as well as the boundedness, compactness, and behavior of their products. These properties exhibit certain parallels with those of classical Toeplitz operators. We begin with a discussion on the uniqueness of symbols.

	\begin{theorem}\label{2Th7}
		The dual truncated Hankel operator \hp ~is a zero operator only when the symbol $\phi$ is identically zero.
	\end{theorem}
	\begin{proof}
		Suppose $ \mathcal{H}_\phi^\theta f = 0 $  for every $f \in K_\theta^\perp.$ Then , $ Q_\theta J(\phi f)= 0 $ implies $ \phi f \in J K_\theta ,$ for any $f\in K_\theta^\perp.$ Choosing $f= \bar{z}^n $ for $n\geq1$, we have $\phi \in z^n JK_\theta$ for every $n\geq 1.$ This implies $J(\phi) = \breve{\phi} \in \bigcap_{n=1}^\infty \bar{z}^n K_\theta = \{0\}. $ And hence $\phi$ is 0. This establishes the uniqueness of the symbol of DTHO. 
	\end{proof}
Next, we turn our attention to the boundedness and compactness of DTHOs.
	\begin{theorem}\label{2Th8}
		The dual truncated Hankel operator \hp ~is bounded if and only if the symbol $\phi$ is an $L^\infty$- element, and  $\|\mathcal{H}_{\phi}^{\theta}\| = ||\phi||_\infty$. Moreover, the only compact DTHO is the zero operator.
	\end{theorem}
	\begin{proof}
		Note that, $\mathcal{H}_\phi h = Q_\theta J (\phi h)= QJ(\phi h) + \theta P \bar{\theta}J(\phi h),$ for $h\in K_\theta^\perp.$ Therefore, for every $h\in K_\theta^\perp$ we get
		$$ ||QJ(\phi h)|| \leq || \mathcal{H}_\phi h || \implies || J QJ(\phi h)|| \leq || \mathcal{H}_\phi h ||,$$ 
		and hence $||T_{\phi\bar{z}} h|| \leq ||\mathcal{H}_\phi h||.$ So, $\mathcal{H}_\phi$ is bounded implies the Toeplitz operator $T_{\bar{z}\phi}$ is bounded. Therefore, $\phi$ is an $L^\infty$ element. And, the converse follows easily.
		\vspace{0.01in}
		
		For compactness, let us consider a weakly convergent sequence $ \{f_n\}$  in $  H^2, $ converging to $0$. Letting $\mathcal{H}_\phi$ to be a compact operator we have $ ||\mathcal{H}_\phi (\theta f_n)|| \rightarrow 0 $ as $n\rightarrow \infty$ implies that $|| T_{\bar{z}\phi\theta} f_n || \rightarrow 0,$ as $n\rightarrow \infty$,  that is the Toeplitz operator $T_{\bar{z}\phi\theta}$ is a compact operator. Therefore, the symbol $\phi$ of DTHO is 0. This completes the proof. 
	\end{proof}
We now examine a property of the product of DTHOs.
	\begin{theorem}\label{2L3}
		Let $\phi_i, \psi_i$ ($i=1,2,\ldots,n$) are $L^\infty$- elements. Then, for any compact operator $K$ on $K_\theta^ \perp$, we have 
		$$  \left\| \sum_{i=1}^n \mathcal{H}_{\phi_i}^\theta \mathcal{H}_{\psi_i}^\theta + K \right\|\geq \left\| \sum_{i=1}^n \breve{\phi_i} \psi_i \right\|_\infty.$$
	\end{theorem}
	\begin{proof}
		Consider unit vectors $f,g \in H^2$. Then for any $m\geq 0$, we have the following.
		\begin{align*}
			& \left \langle \left(\sum_{i=1}^n \mathcal{H}_{\phi_i} \mathcal{H}_{\psi_i} + K\right) (\theta z^m f), \theta z^m g \right \rangle \\
			& = \sum_{i=1}^n\langle \mathcal{H}_{\psi_i} (\theta z^m f), \mathcal{H}_{\phi_i}^*(\theta z^m g) \rangle + \langle K(\theta z^m f ), \theta z^m g \rangle\\
			& = \sum_{i=1}^n \langle Q J (\psi_i \theta z^m f), Q J (\phi_i^* \theta z^m g) \rangle + \sum_{i=1}^n \langle \theta P \bar{\theta} J(\psi_i z^m f ), \theta P \bar{\theta} J(\phi_i^* z^m g) \rangle \\
			& \hspace{0.5in} + \langle K(\theta z^m f ), \theta z^m g \rangle\\
			& \longrightarrow \sum_{i=1}^n \langle \psi_i f, \phi_i^* g \rangle + 0 + 0 \text{ as } m \rightarrow \infty
		\end{align*}
		because $P( \bar{\theta} J (\psi \theta z^m  f )) \rightarrow 0$, and $K(\theta z^m f ) \rightarrow 0 $ since $K$ is a compact operator on $K_\theta^\perp $, and $\langle Q J (\psi_i \theta z^m f), Q J (\phi_i^* \theta z^m g) \rangle \rightarrow \langle \psi_i f, \phi_i^* g \rangle $ as $m \rightarrow\infty.$  Therefore, $||\sum_{i=1}^n \mathcal{H}_{\phi_i} \mathcal{H}_{\psi_i} + K|| \geq |\sum_{i=1}^n \langle \psi_i f, \phi_i^* g \rangle | = |\sum_{i=1}^n \langle \breve{\phi_i} \psi_i f, g \rangle |=| \sum_{i=1}^n \langle T_{\breve{\phi_i} \psi_i} f, g \rangle| = \left|\left\langle T_{\sum_{i=1}^n\breve{\phi_i} \psi_i}(f), g\right\rangle \right|.$  Taking the supremum over the unit vectors $ f$, $g$, we have the desired result.
	\end{proof}
	As a consequence of the above theorem, we have the following result.
	\begin{corollary}\label{2C7}
		If the product $\mathcal{H}_\phi^\theta \mathcal{H}_\psi ^\theta $ is compact, then $ \breve{\phi} \psi = 0. $
	\end{corollary}

    In particular, the product  \hp$\mathcal{H}_\psi^\theta=0$ implies that $\breve{\phi}\psi=0.$ Now, it is natural to inquire whether further conclusions can be established for each individual DTHO. Proceeding with the elegant method by Ding and Sang \cite{DS}, and using a relation between $T_\phi$ and $B_\psi$ (see \eqref{Toep}), we deduce the following conclusion. 

\begin{theorem}
    The product of two dual truncated Hankel operators \hp ~and $ \mathcal{H}_\psi ^\theta $ is zero if and only if one of them is zero.
\end{theorem}

\begin{proof}
    It is clear that if the symbol of one of the DTHO is zero, then their product is zero and hence $\breve{\phi}
    \psi=0$. 

    For the converse, let $\h_\phi \h_\psi=0$, then we have $ \breve{\phi}\psi=0. $
Suppose $\phi|_E=0$, for a non-zero Lebesgue measurable subset $E$ of $\mathbb{T}$. Then $\h_\phi \h_\psi=0$ implies $$(I-P)JM_\phi(\h_\psi h)+ M_\theta P M_{\overline{\theta}} J M_\phi (\h_\psi h)=0$$ for every $h\in K_\theta^\perp.$
Therefore, \begin{align*}
    & (I-P)JM_\phi(\h_\psi h)=0 \text{ for every } h\in K_\theta^\perp \\
    & \implies (I-P)JM_\phi\big( (I-P)JM_\psi h + M_\theta P M_{\overline{\theta}}JM_\psi h \big)=0\\
    & \implies ( JM_\phi - PJM_\phi ) \big( JM_\psi h - PJM_\psi h + M_\theta P M_{\overline{\theta}}JM_\psi h \big)=0  
    \end{align*}
    \begin{align*}
    & \implies JM_\phi JM_\psi h - PJM_\phi JM_\psi h - JM_\phi PJM_\psi h + PJM_\phi PJM_\psi h + \\ 
    & \hspace{1in}JM_\phi  M_\theta P M_{\overline{\theta}}JM_\psi h - PJM_\phi M_\theta P M_{\overline{\theta}}JM_\psi h = 0\\
    & \implies M_{\breve{\phi}} M_\psi h - P M_{\breve{\phi}} M_\psi h - JM_\phi PJM_\psi h + PJM_\phi PJM_\psi h + \\ 
    & \hspace{0.5in}JM_\phi  M_\theta P M_{\overline{\theta}}JM_\psi h - PJM_\phi M_\theta P M_{\overline{\theta}}JM_\psi h = 0 \text{ (since } JM_\phi J= M_{\breve{\phi}})\\
    &\implies JM_\phi PJM_\psi h - JM_\phi  M_\theta P M_{\overline{\theta}}JM_\psi h  = PJM_\phi PJM_\psi h - \\
    &\hspace{2in} PJM_\phi M_\theta P M_{\overline{\theta}}JM_\psi h, \text{ for every } h\in K_\theta^\perp.
\end{align*}

Now consider $h=\breve{\theta}f$ for any $f\in \overline{H_0^2},$ and, hence, $h\in \overline{H_0^2}$. Thus, we have
\begin{align*}
    & JM_\phi PJM_\psi \breve{\theta}f - JM_\phi  M_\theta P M_{\overline{\theta}}JM_\psi \breve{\theta}f  = PJM_\phi PJM_\psi \breve{\theta}f \\
    & \hspace{3in}-PJM_\phi M_\theta P M_{\overline{\theta}}JM_\psi \breve{\theta}f\\
    & \implies JM_\phi (B_{\breve{\theta}\psi}-T_\theta B_{\psi})(f)= (H_\phi B_{\breve{\theta}\psi} - H_{\theta \phi} B_\psi )(f)\\
    & \implies (H_\phi B_{\breve{\theta}\psi} - H_{\theta \phi} B_\psi )(f)=0, 
    \end{align*}
because the right-hand side of the above equality is an $H^2$ element and the left-hand side equals $0$ on a non-zero Lebesgue measurable subset of $\mathbb{T}$. This further implies
    $$ H_\phi ( B_{\breve{\theta}\psi}- T_\theta B_\psi )(f) )= 0.$$
On the other hand, it is not difficult to see that $ker{H_\phi}=\{0\}.$ Therefore, we have 
$$ ( B_{\breve{\theta}\psi}- T_\theta B_\psi )(f) )= 0, \text{ for any $f\in \overline{H_0^2}$ } .$$ So, we conclude that $\psi\in z\overline{H^2}$ (because $T_\phi B_\psi = B_{\breve{\phi}\psi}$ if and only if either $\phi$ is co-analytic or $\psi \in z\overline{H^2}$). And since, $\breve{\phi}\psi=0$ and $\phi|_{E}=0,$ we get $\psi=0.$ This completes the proof.
\end{proof}
	
	\subsection{Hankelness of DTHO}
    
In the previous subsection, we observed that DTHOs share certain algebraic properties with classical Toeplitz operators, since some of their components behave like Toeplitz operators. At the same time, other components exhibit behavior similar to Hankel operators. This subsection addresses the question posed in \eqref{2Q} in the introduction and provides an affirmative answer under certain restrictions. \emph{Throughout this subsection, we consider the inner function $\theta$ is symmetric, that is, $\theta^*(z)= \overline{\theta(\bar{z})} =\theta(z).$}

	\begin{lemma}\label{2L4}
		For two dual truncated Hankel operators $\mathcal{H}_\phi ^\theta$ and $\mathcal{H}_\psi ^\theta$, where $\phi,\psi \in Z_\theta$ and $\theta_0=0$, the following identity holds.
		$$ \mathcal{H}_\phi^\theta \mathcal{H}_\psi^\theta - U_\theta^*\mathcal{H}_\phi^\theta \mathcal{H}_\psi^\theta U_\theta = Q_\theta J(\phi\theta) \otimes Q_\theta(\psi^* \theta) .$$
	\end{lemma}
	
	\begin{proof}
		Combining Lemma~\ref{2L1}, Theorem~\ref{2Th1}, and Corollary~\ref{2C1} we get
		\begin{align*}
			& \mathcal{H}_\phi \mathcal{H}_\psi - U^* \mathcal{H}_\phi \mathcal{H}_\psi U \\
			& = \mathcal{H}_\phi \mathcal{H}_\psi - \mathcal{H}_\phi U U^* \mathcal{H}_\psi\\
			& = \mathcal{H}_\phi (I- UU^*) \mathcal{H}_\psi = \mathcal{H}_\phi (\theta \otimes \theta) \mathcal{H}_\psi\\
			& = \mathcal{H}_\phi(\theta) \otimes \mathcal{H}_\psi^*(\theta) = Q_\theta J(\phi\theta) \otimes Q_\theta J(\psi^* \theta) .
		\end{align*}
		This completes the proof. 
	\end{proof}
	We know that the product of two classical Hankel operators satisfies the Brown-Halmos  \cite{BH} Toeplitz identity if and only if one of the Hankel operators is zero \cite{AR}. The following theorem establishes an analog of this fact in the context of DTHO with some restrictions. 
	\begin{theorem}
		The operator $\mathcal{H}_\phi ^\theta \mathcal{H}_\psi ^\theta$ with the condition that $\phi,\psi \in Z_\theta$ and $\theta_0=0$, satisfies the Brown-Halmos Toeplitz identity on $K_\theta^\perp$ if and only if either $\phi=0$ or $\psi=0.$
	\end{theorem}

	\begin{proof}
    
		By the Brown--Halmos Toeplitz identity, an operator $T$ on $H^2$ is Toeplitz if and only if $S^*TS=T$. So, in our context,  the product $\mathcal{H}_\phi  \mathcal{H}_\psi $ satisfies the Brown-Halmos Toeplitz identity on $K_\theta^\perp$ if and only if $$ \mathcal{H}_\phi  \mathcal{H}_\psi = U^* \mathcal{H}_\phi  \mathcal{H}_\psi U .$$
		For $\phi,\psi \in Z_\theta, \text{ and } \theta_0=0,$ by Lemma~\ref{2L4}, we conclude that
		\begin{align*}
			&	\mathcal{H}_\phi \mathcal{H}_\psi - U^* \mathcal{H}_\phi \mathcal{H}_\psi U= 0  \iff Q_\theta J(\phi\theta) \otimes Q_\theta J(\psi^* \theta) = 0\\
			& \iff \text{ either }Q_\theta J(\phi\theta) =0 \text{ or }  Q_\theta(\psi^* \theta) = 0 \\
			& \iff \text{ either } \phi \theta \in J K_\theta = \overline{K_\theta} \text{ or }  \psi^* \theta \in J K_\theta = \overline{K_\theta} \text{ ( since, $JK_\theta = \overline{K_\theta}$ as $ \theta^*=\theta $)}.
		\end{align*}
		On the other hand, $ \phi = cz + \bar{h}_1 \implies \theta \phi = cz\theta + \theta \bar{h}_1 \in \overline{K_\theta}.$ And, for any $h\in \overline{K_\theta}$ $\langle \theta\phi , \bar{h} \rangle =\langle cz\theta + \theta \bar{h}_1 , \bar{h} \rangle = \langle \theta \bar{h}_1, \bar{h} \rangle = \langle \theta h, h_1 \rangle=0.$ This implies $\phi\theta\in \overline{K_\theta}^{\perp}\cap \overline{K_\theta}$, which implies $\phi =0$ since $\theta \neq 0 $ a.e. on $\mathbb{T}.$ Similarly, one can also show by considering the other condition that $ \psi = 0$.
	\end{proof}
	Next, we discuss the commutativity of DTHOs under the same conditions as above.

	\begin{theorem}
		Let $\phi,\psi \in Z_\theta$ and $\theta_0=0$. Then $\mathcal{H}_\phi ^\theta$ and $\mathcal{H}_\psi ^\theta$   commute if and only if one of them is a scalar multiple of the other.
	\end{theorem}
	
	\begin{proof}
		According to Lemma \ref{2L4} we have $$\mathcal{H}_\phi \mathcal{H}_\psi - U^* \mathcal{H}_\phi \mathcal{H}_\psi U = Q_\theta J(\phi\theta) \otimes Q_\theta J(\psi^* \theta), $$ and 		
		$$\mathcal{H}_\psi \mathcal{H}_\phi - U^* \mathcal{H}_\psi \mathcal{H}_\phi U = Q_\theta J(\psi\theta) \otimes Q_\theta J(\phi^* \theta) .$$
		Therefore, $\mathcal{H}_\phi \mathcal{H}_\psi = \mathcal{H}_\psi \mathcal{H}_\phi$ implies that \begin{align}
		    & Q_\theta J(\phi\theta) \otimes Q_\theta J(\psi^* \theta) = Q_\theta J(\psi\theta) \otimes Q_\theta J(\phi^* \theta),
		\end{align}
		and hence there exists a complex number $c$ such that $ Q_\theta J(\phi\theta) = c Q_\theta J(\psi\theta) $ and $ \bar{c} Q_\theta J(\psi^*\theta) =  Q_\theta J(\phi^* \theta).$ Now considering $ Q_\theta J(\phi\theta) = c Q_\theta J(\psi\theta)  \implies  Q_\theta J((\phi - c\psi)\theta) = 0. $ Therefore, by an argument analogous to that of the previous theorem, we obtain $\phi=c\psi,$ and hence the corollary follows. 
	\end{proof}
	We get the following result as a direct consequence of the above theorem.
	\begin{corollary}\label{2normal}
		Every normal dual truncated Hankel operator $\mathcal{H}_\phi^\theta$  with  $\phi\in Z_\theta$ and $\theta_0=0$ is merely some scalar multiple of the self-adjoint ones.    
	\end{corollary}
	\begin{proof}
		Let $\mathcal{H}_\phi$ be a normal operator. Then $ \mathcal{H}_\phi^*\mathcal{H}_\phi= \mathcal{H}_\phi^* \mathcal{H}_\phi $, which by using the previous theorem implies $ \mathcal{H}_\phi = c \mathcal{H}_{\phi^*}$ for some complex number $c.$ Since $|| \mathcal{H}_\phi ||=|| \mathcal{H}_{\phi^*} ||, $ the complex number $ c$ is unimodular. Taking $c= e^{i\xi}$, leads to $( e^{\frac{-i}{2}\xi} \mathcal{H}_\phi )^* = e^{\frac{-i}{2}\xi} \mathcal{H}_\phi, $ hence the result follows. 
	\end{proof}
Next, we turn our attention to the hyponormality of DTHOs. Before proceeding, we briefly recall the definition of a hyponormal operator. An operator $T: H\rightarrow H,$ where $H$ is a Hilbert space, is called hyponormal if  $ T^*T-TT^*$ is a positive operator. In \cite{CC}, Cowen provided the first characterization of hyponormality for classical Toeplitz operators. Since then, several mathematicians have extended the study of hyponormality to various classes of operators \cite{CHL, CLY, CG4, CG3, CG2}.

It is well known that classical hyponormal Hankel operators on the Hardy space $H^2$ are, in fact, normal. In the following theorem, we establish a similar property for DTHOs under the condition $\theta^*=\theta,$ regardless of the choice of symbol.

	\begin{theorem}\label{2hypon}
		Every hyponormal dual truncated Hankel operator $\mathcal{H}_\phi ^\theta$ is normal. 
	\end{theorem}
	
	\begin{proof}
		In order to prove this theorem, we need the relation $ || \mathcal{H}_\phi^* f^* ||=|| \mathcal{H}_\phi f || $ for every $f \in K_\theta^\perp.$ Now, from the definition of DTHO, it follows
		\begin{align*}
			|| \mathcal{H}_\phi^* f^* ||^2 & = \langle \mathcal{H}_\phi^* f^* , \mathcal{H}_\phi^* f^* \rangle  \\
			& = \langle Q_\theta J (\phi^* f^*) , Q_\theta J (\phi^* f^*)  \rangle\\
			& =\langle  Q_\theta J (\phi^* f^*) ,J( \phi^* f^*)  \rangle \\
			& = \langle Q_\theta (\overline{\phi f}) , \overline{ \phi f} \rangle \\
			&= \langle  \phi f, \overline{Q_\theta (\overline{\phi f}} ) \rangle \\
			&= \langle \phi f, \overline{ Q(\overline{\phi f})} + \overline{\theta P \bar{\theta} (\overline{\phi f}) } \rangle\\
			& =\langle \phi f, \overline{ Q(\overline{\phi f})} + \langle \phi f, \bar{\theta}\overline{P (\overline{\phi \theta f}) } \rangle\\
			& = \langle J (\phi f), QJ(\phi f) \rangle + \langle \theta \phi f,\overline{P (\overline{\phi \theta f}) } \rangle\\
			& =\langle J(\phi f), QJ(\phi f) \rangle + \langle J( \theta \phi f) ,  PJ(\theta \phi f)\rangle\\
			&  =\langle J(\phi f), QJ(\phi f) \rangle + \langle J( \phi f) ,  \theta^* P\breve{\theta }J(\phi f)\rangle\\
			& =\langle J(\phi f), QJ(\phi f) \rangle + \langle J( \phi f) ,  \theta P\bar{\theta }J(\phi f)\rangle\\
			& =\langle  J(\phi f) , Q_\theta J(\phi f) \rangle  = \langle  Q_\theta(J(\phi f)) , Q_\theta J(\phi f) \rangle\\
			& = \langle \mathcal{H}_\phi f,\mathcal{H}_\phi f\rangle  =||\mathcal{H}_\phi f||^2.
		\end{align*}
		Therefore, $ || \mathcal{H}_\phi^* f^* ||=|| \mathcal{H}_\phi f || $ for every $f \in K_\theta^\perp.$ Now, replacing $ g^*$ in place of $f,$ we have $ || \mathcal{H}_\phi^* g ||=|| \mathcal{H}_\phi g^* || .$ Suppose that $\mathcal{H}_\phi $ is hyponormal. Then $ || \mathcal{H}_\phi f || \geq || \mathcal{H}_\phi^* f || $ holds for every $f\in K_\theta^\perp$. This implies $ ||\mathcal{H}_\phi g^* || \geq || \mathcal{H}_\phi^* g^* || \implies ||\mathcal{H}_\phi^* g|| \geq ||\mathcal{H}_\phi g||.$ On the other hand, by the hyponormality of $\mathcal{H}_\phi$ we have $ || \mathcal{H}_\phi g || \geq || \mathcal{H}_\phi^* g ||.$ So, $ || \mathcal{H}_\phi f ||=||  \mathcal{H}_\phi^* f ||$ for every $f \in K_\theta^\perp.$ Hence, \hp ~is normal. This completes the proof.

	\end{proof}
\textbf{{\emph{Special Remark:}}}
\vspace{0.1in}

(1) In Theorem~\ref{2hypon}, we observed that hyponormal DTHOs (that is, those with $\theta^*=\theta$) are, in fact, normal. In Corollary~\ref{2normal}, we established the normality of DTHOs under certain restrictions. It is, therefore, natural and interesting to investigate the conditions under which a DTHO is normal without any such restrictions.
\vspace{0.1in}

(2)Another interesting direction of investigation is determining the conditions under which the product of two DTHOs is itself a DTHO or a DTTO.

	 \section*{\textbf{{Acknowledgment}}}
	We sincerely thank the anonymous referee for the careful reading of the manuscript and for the valuable suggestions. A. Chattopadhyay is supported by the Core Research Grant (CRG), File No: CRG/2023/004826, by the Science and Engineering Research Board(SERB), Department of Science \& Technology (DST), and the FIST Grant No.: SR/FST/MS-II/2024/183 (C) by the Department of Science \& Technology (DST), Government of India. S. Jana gratefully acknowledges the support provided by IIT Guwahati, Government of India.

    \noindent $^*$ [A. Chattopadhyay] Department of Mathematics, Indian Institute of Technology Guwahati, Guwahati, 781039, India.\\
		\textit{Email address:}~ arupchatt@iitg.ac.in, 2003arupchattopadhyay@gmail.com.
		\vspace{0.05in}
	
		\noindent $^{**}$ [S. Jana] Department of Mathematics, Indian Institute of Technology Guwahati, Guwahati, 781039, India.\\
		\textit{Email address:}~ supratimjana@iitg.ac.in, suprjan.math@gmail.com.
	
\end{document}